\documentstyle[12pt]{article}
\begin{document}
\newcommand{\text}[1]{\mbox{{\rm #1}}}  
\newcommand{\gd}{\delta}
\newcommand{\itms}[1]{\item[[#1]]} 
\newcommand{\nin}{\in\!\!\!\!\!/}
\newcommand{\sub}{\subset} 
\newcommand{\cntd}{\subseteq}     
\newcommand{\go}{\omega} 
\newcommand{\Pa}{P_{a^\nu,1}(U)} 
\newcommand{\fx}{f(x)}  
\newcommand{\fy}{f(y)} 
\newcommand{\gD}{\Delta}
\newcommand{\gl}{\lambda} 
\newcommand{\gL}{\Lambda} 
\newcommand{\half}{\frac{1}{2}} 
\newcommand{\sto}[1]{#1^{(1)}}
\newcommand{\stt}[1]{#1^{(2)}}
\newcommand{\Z}{\hbox{\sf Z\kern-0.720em\hbox{ Z}}}
\newcommand{\singcolb}[2]{\left(\begin{array}{c}#1\\#2
\end{array}\right)} 
\newcommand{\ga}{\alpha}
\newcommand{\gb}{\beta} 
\newcommand{\gga}{\gamma}
\newcommand{\ul}{\underline} 
\newcommand{\ol}{\overline} 
\newcommand{\qed}{\kern 5pt\vrule height8pt width6.5pt depth2pt}
\newcommand{\Lrraro}{\Longrightarrow}
\newcommand{\Nb}{|\!\!/}
\newcommand{\NN}{{\rm I\!N}}
\newcommand{\bsl}{\backslash}     
\newcommand{\gt}{\theta}
\newcommand{\op}{\oplus}
\newcommand{\Op}{\bigoplus}          
\newcommand{\CR}{{\cal R}}
\newcommand{\tr}{\bigtriangleup}
\newcommand{\grr}{\omega_1} 
\newcommand{\ben}{\begin{enumerate}}
\newcommand{\een}{\end{enumerate}}
\newcommand{\ndiv}{\not\mid}
\newcommand{\bab}{\bowtie}
\newcommand{\hal}{\leftharpoonup}
\newcommand{\har}{\rightharpoonup}
\newcommand{\ot}{\otimes}
\newcommand{\OT}{\bigotimes}
\newcommand{\bwe}{\bigwedge}
\newcommand{\gep}{\varepsilon}
\newcommand{\gs}{\sigma} 
\newcommand{\rbraces}[1]{\left( #1 \right)}
\newcommand{\bbox}{$\;\;\rule{2mm}{2mm}$}
\newcommand{\sbraces}[1]{\left[ #1 \right]}
\newcommand{\bbraces}[1]{\left\{ #1 \right\}}
\newcommand{\OO}{_{(1)}}
\newcommand{\TT}{_{(2)}}
\newcommand{\FF}{_{(3)}}
\newcommand{\minus}{^{-1}}
\newcommand{\CV}{\cal V} 
\newcommand{\CVs}{\cal{V}_s} 
\newcommand{\un}{U_q(sl_n)'}
\newcommand{\on}{O_q(SL_n)'}
\newcommand{\slq}{U_q(sl_2)}
\newcommand{\olq}{O_q(SL_2)}
\newcommand{\UU}{U_{(N,\nu,\go)}}
\newcommand{\HH}{H_{n,q,N,\nu}} 
\newcommand{\ct}{\centerline}
\newcommand{\bs}{\bigskip}
\newcommand{\qua}{\rm quasitriangular}   
\newcommand{\ms}{\medskip}
\newcommand{\noin}{\noindent}
\newcommand{\mat}[1]{$\;{#1}\;$}
\newcommand{\raro}{\rightarrow}
\newcommand{\map}[3]{{#1}\::\:{#2}\raro{#3}}
\newcommand{\alg}{{\rm Alg}}
\def\newtheorems{\newtheorem{theorem}{Theorem}[subsection]
                 \newtheorem{cor}[theorem]{Corollary}
                 \newtheorem{prop}[theorem]{Proposition}
                 \newtheorem{lemma}[theorem]{Lemma}
                 \newtheorem{defn}[theorem]{Definition}
                 \newtheorem{Theorem}{Theorem}[section]
                 \newtheorem{Corollary}[Theorem]{Corollary}
                 \newtheorem{Proposition}[Theorem]{Proposition}
                 \newtheorem{Lemma}[Theorem]{Lemma}
                 \newtheorem{Defn}[Theorem]{Definition}
                 \newtheorem{Example}[Theorem]{Example}
                 \newtheorem{Remark}[Theorem]{Remark} 
                 \newtheorem{claim}[theorem]{Claim}
                 \newtheorem{sublemma}[theorem]{Sublemma}
                 \newtheorem{example}[theorem]{Example}
                 \newtheorem{remark}[theorem]{Remark}
                 \newtheorem{question}[theorem]{Question}
                 \newtheorem{conjecture}{Conjecture}[subsection]}
\newtheorems
\newcommand{\proof}{\par\noindent{\bf Proof:}\quad}
\newcommand{\dmatr}[2]{\left(\begin{array}{c}{#1}\\
                            {#2}\end{array}\right)}
\newcommand{\doubcolb}[4]{\left(\begin{array}{cc}#1&#2\\
#3&#4\end{array}\right)}
\newcommand{\qmatrl}[4]{\left(\begin{array}{ll}{#1}&{#2}\\
                            {#3}&{#4}\end{array}\right)}
\newcommand{\qmatrc}[4]{\left(\begin{array}{cc}{#1}&{#2}\\
                            {#3}&{#4}\end{array}\right)}
\newcommand{\qmatrr}[4]{\left(\begin{array}{rr}{#1}&{#2}\\
                            {#3}&{#4}\end{array}\right)}
\newcommand{\smatr}[2]{\left(\begin{array}{c}{#1}\\
                            \vdots\\{#2}\end{array}\right)}

\newcommand{\ddet}[2]{\left[\begin{array}{c}{#1}\\
                           {#2}\end{array}\right]}
\newcommand{\qdetl}[4]{\left[\begin{array}{ll}{#1}&{#2}\\
                           {#3}&{#4}\end{array}\right]}
\newcommand{\qdetc}[4]{\left[\begin{array}{cc}{#1}&{#2}\\
                           {#3}&{#4}\end{array}\right]}
\newcommand{\qdetr}[4]{\left[\begin{array}{rr}{#1}&{#2}\\
                           {#3}&{#4}\end{array}\right]}

\newcommand{\qbracl}[4]{\left\{\begin{array}{ll}{#1}&{#2}\\
                           {#3}&{#4}\end{array}\right.}
\newcommand{\qbracr}[4]{\left.\begin{array}{ll}{#1}&{#2}\\
                           {#3}&{#4}\end{array}\right\}}

\title{On Semisimple Hopf Algebras of Dimension $pq$}
\author{Shlomo Gelaki $^{1}$ 
\\Department of Mathematics and Computer Science\\
Ben Gurion University of  the Negev\\Beer Sheva, Israel
\and Sara Westreich $^1$\\Interdisciplinary Department of the Social Science\\
Bar-Ilan University\\Ramat-Gan, Israel}
\footnotetext[1]{This research was supported by THE ISRAEL SCIENCE FOUNDATION
founded by the Israel Academy of Sciences and Humanities.}
\date{July 23, 1997}
\maketitle

In this paper 
we consider semisimple Hopf algebras of 
dimension $pq$ over an algebraically closed field $k$ of characteristic $0,$
where $p$ and $q$ are distinct prime numbers. 
Masuoka has proved that a semisimple Hopf algebra of dimension $2p$ over
$k,$
where $p$ is an odd prime, is either a group algebra or a dual of a group 
algebra [Ma1]. The authors have pushed the analysis further and obtained 
the same
result for semisimple $A$ of dimension $3p,$ where $p$ is a prime greater 
than $3$ [GW]. Thus, a natural conjecture is: \\
{\bf Conjecture 1:} 
Any semisimple Hopf algebra of 
dimension $pq$ over $k,$ where $p$ and $q$ are distinct prime numbers, is 
either a group algebra or a dual of a group algebra.

A well known property of $A,$ a finite dimensional semisimple group
algebra or a dual of a group algebra, is that it is of {\em Frobenius type};
that is,
the dimension of any irreducible representation of
$A$ divides the dimension of $A$ 
(the definition is due to Montgomery [Mo2]).
A special case of one of Kaplansky's conjectures [Kap] is:\\
{\bf Conjecture 2:}
Any semisimple Hopf algebra of 
dimension $pq$ over $k,$ where $p$ and $q$ are distinct prime numbers, is 
of Frobenius type.\\
In this paper we prove among the rest that
Conjecture 1 is equivalent to Conjecture 2 ({\bf Theorem \ref{rem}}).

A major role in the analysis is played by $G(A),$ the group of grouplike
elements of $A.$
By [NZ], $|G(A)|$ is either $1,p,q$ or $pq.$
We prove in {\bf Theorem \ref{notq}} that if $p<q$ then $|G(A)|\ne q,$ and 
if $|G(A)|=p$ then $q=1(mod\,p).$ Consequently, we prove in {\bf Theorem
\ref{trivial}} that if $|G(A)|\ne 1$ and $q\ne 1(mod\,p)$ then $A$ is a 
group algebra. 

Thus Theorem \ref{trivial} suggests a question: When is
$|G(A)|\ne 1$? In
{\bf Proposition \ref{nt}} we prove that this is guarenteed when
$A^*$ is of Frobenius type, 
and in {\bf Theorem \ref{frob}} we prove that if moreover $q\ne
1(mod\,p)$ then $A$ is a group algebra. 
In {\bf Theorem \ref{main}} we prove that if
$|G(A^*)|\ne 1$ and $A^*$ is of Frobenius type, then $A$ is
either a group algebra or the dual of a group algebra, and
$|G(A)|=p<q$ or $pq.$ 
The equivalence of Conjectures 1 and 2 is thus a consequence of Proposition
\ref{nt} and Theorem \ref{main}.

We conclude by using Theorem \ref{trivial} to prove
in {\bf Theorem \ref{5}} that if $A$ is a semisimple Hopf algebra
of dimension $5p,$ $p$ an odd prime, and if 
$p=2$ or $4(mod\,5)$ or $p\in\{13,23\}$, then $A$
is a group algebra. Moreover, we obtain in {\bf Theorem \ref{7}} the
same result for semisimple $A$ of dimension $7p,$ $p$ a prime, and
$p=6(mod\,7)$ or $p\in\{17,31\}.$ 
\section{Preliminaries}
In this paper $k$ will always denote
an algebraically closed field of characteristic $0.$

Recall that a finite dimensional Hopf algebra over $k$ is semisimple if
and only if it is cosemisimple [LR].

Let $A$ be a finite dimensional semisimple Hopf algebra over $k,$ and
let $\rho_{_V}:A\raro End_k(V)$ be a finite dimensional
representation of $A.$ The afforded character $\chi_{_V}$ 
is given by $\chi _{_V}(a)=tr(\rho_{_V}(a))$ for all $a\in A.$ 
A character $\chi _{_V}$ is called irreducible if the representation $V$
is irreducible.
Let $C(A)$ denote the character ring of $A;$ that is, the $k$-subalgebra
of $A^*$ generated by the characters $\chi_{_V}$ of finite
dimensional $A$-modules $V.$
The set of all the irreducible characters forms a basis of $C(A)$ [La]. 
Zhu has proved that $C(A)$ is semisimple 
and if $e_{A^*},e_1,\dots ,e_k$ are the primitive idempotents of $C(A),$ 
where $e_{A^*}$ is an integral of $A^*,$ then 
\begin{equation}\label{zhu}
dim\,A=1+\sum_{i=1}^k dim(e_iA^*) 
\end{equation}
and the dimension of each
$e_iA^*$ divides the dimension of $A$ [Z]. Note that $dimC(A)\ge k+1,$
and equality holds if and only if $C(A)$ is commutative.\\
Let $f:A\raro A^*$ be the map given by $f(a)=a\har \gl=\sum <a,
\gl _{(2)}>\gl _{(1)}$ for all $a\in A,$ where $\gl$ is a non-zero 
integral of
$A^*.$ Recall that $f$ gives a linear isomorphism  
between $kG(A)$ and the sum of the $1-$dimensional ideals of $A^*$ 
(where $G(A)$ denotes the group of grouplike elements of $A$),
and a linear
isomorphism between the center $Z(A)$ of $A$ and $C(A).$ Therefore, using
the notation of (\ref{zhu}), $dim(e_iA^*)=1$ for some $i$ if and only if
$G(A)\cap Z(A)\ne \{1\}.$

Let $A$ be a finite dimensional semisimple Hopf algebra over $k.$
Any simple subcoalgebra $C_l$ of $A^*$ has a basis
$\{x_{ij}^l|1 \le i,j \le n_l\},$ where $\gD(x_{ij}^l) = \sum_{k=1}^{n_l}
x_{ik}^l \ot x_{kj}^l$ and $\gep(x_{ij}^l)=\gd_{i,j}.$ 
Note that $n_l=1$ if and only if $C_l=\{g\}$ for some $g\in G(A).$ 
Nichols and Richmond have proved that if $dimA$ is odd then  
$A$ does not have a $2-$dimensional irreducible module [NR], hence
\begin{equation}\label{square}
dimA=|G(A)|+\sum_ln_l^2,\qquad n_l\ge 3. 
\end{equation}
Now, $L$ is an irreducible left
coideal of $C_l$ if and only if  
\begin{equation}\label{left}
L=L_j^l=sp\{x_{kj}^l|1 \le k \le n_l\}
\end{equation}
for some $1 \le j \le n_l.$
Similarly, 
$R$ is an irreducible right
coideal of $C_l$ if and only if  
\begin{equation}\label{right}
R=R_k^l=sp\{x_{kj}^l|1 \le j \le n_l\}
\end{equation}
for some $1 \le k \le n_l.$
Note that 
\begin{equation}\label{inter}
dim(L_j^l\cap R_k^l)=1 
\end{equation}
for any $1 \le j,k \le n_l.$ 

The following theorem will be very useful in the sequel.
\begin{Theorem} {\bf [R]}\label{radf}
If $H \stackrel{i}{\hookrightarrow} A \stackrel{\pi}
{\raro} H$ is a sequence of finite dimensional 
Hopf algebra maps where $i$ is injective,
$\pi$ is surjective and $\pi \circ i=id_H$ then there exists $B\subset A$
so that: 
\begin{description}
\item [(i)] $B$ is a left $H-$module algebra and coalgebra via the adjoint
action.
\item [(ii)] $B$ is a left $H-$comodule algebra and coalgebra via $\rho (b)=
\sum b^{(1)}\ot b^{(2)}=\sum\pi(b_{(1)})\ot b_{(2)}.$
\item [(iii)] $B\cong A/AH^{+}$ as a coalgebra, via the map $b\times h
\mapsto b\gep (h).$
\item [(iv)] $B$ is a subalgebra coideal of $A.$
\item [(v)]  As an algebra $A=B \times H$ is a smash product.
\item [(vi)] As a coalgebra $A=B \times H$ is a smash coproduct,
that is: $\gD(b\times h)=\sum b_{(1)}\times b_{(2)}^{(1)}h_{(1)}\ot 
b_{(2)}^{(2)}\times h_{(2)}.$
\item [(vii)] The map $B\times H\raro A$ $(b\times h\mapsto bi(h))$ is an
isomorphism of bialgebras. 
\end{description}
\end{Theorem}
\section{Main Results}
We start by determining the order of $G(A).$
\begin{Theorem}\label{notq}
Let $A$ be a semisimple Hopf algebra 
of dimension $pq,$ where $p<q$ are two prime numbers. Then:
\ben
\item $|G(A)|\ne q.$ 
\item If $|G(A)|=p$ then $q=1(mod\,p).$
\een 
\end{Theorem}
\proof 1. Suppose on the contrary that $|G(A)|=q.$ 
If $G(A)\cap Z(A)=G(A)$ then
$H=kG(A)$ is central in $A,$ hence is a 
normal sub Hopf algebra of $A.$ Since $A/(AH^+)$ is a Hopf algebra of 
dimension $p$ it follows by [Z] that $A/(AH^+)\cong kC_p,$ and hence 
$A\cong kC_q\#_{\sigma}kC_p.$ Since $kC_q$ 
is central the afforded action is trivial [Ma3, Section 1]. Therefore $A$ is
isomorphic as an algebra to the twisted group ring $kC_q^t[C_p]$ of the
cyclic group $C_p$ over the commutative algebra $kC_q.$ Hence $A$ must be
commutative (See the proof of [Ma2, Proposition 2.3(2)]). 
Thus, $A^*$ is a group algebra and hence of Frobenius type. By 
(\ref{square}), $pq=q+ap^2+bq^2$ for some integers $a,b\ge 0.$ But $p<q,$ 
hence $b=0,$ which yields a 
contradiction. We conclude that $G(A)\cap Z(A)=\{1\}.$
Therefore, using the notation of (\ref{zhu}), $dim(e_iA^*)\in \{p,q\}$ for 
all $i.$ Let $E_0$ be the integral of $kG(A)$ with $\gep (E_0)=1.$ Since 
$dim(A/AH^+)=p,$ it follows that $AH^+=A(1-E_0)$ has dimension $(q-1)p$ 
and thus $dim(AE_0)=p.$ Moreover, $E_0e_{_A}=e_{_A},$ hence 
$E_0=e_{_A}+\sum_j e_{i_j}.$ But $p<q,$ hence counting
dimensions yields a contradiction and the result follows.\\
2. If $G(A)\cap Z(A)=G(A)$ then
$H=kG(A)$ is central in $A.$ As in the proof of part 1, it follows that 
$A\cong kC_p\#_{\sigma}kC_q,$ and hence that $A$ is
commutative. 
Therefore, $A^*$ is a group algebra and hence of Frobenius type. By 
(\ref{square}), $pq=p+ap^2+bq^2$ for some integers $a,b\ge 0.$ Clearly, 
$b=0$ and hence $q=1+ap.$ If $G(A)\cap Z(A)=\{1\}$
then, using the notation of (\ref{zhu}), $dim(e_iA^*)\in \{p,q\}$ for 
all $i.$ Since 
$dim(A/AH^+)=q,$ it follows that $AH^+=A(1-E_0)$ has dimension $(p-1)q$ 
and thus $dim(AE_0)=q.$ Hence, 
$E_0=e_{_A}+\sum_j e_{i_j}.$ But, counting
dimensions yields that $dim(e_{i_j}A^*)=p$ for all $j,$ and the
result follows in this case as well.\qed\\ 

The following is a direct consequence of Theorem \ref{notq}:
\begin{Theorem}\label{trivial}
Let $A$ be a semisimple Hopf algebra 
of dimension $pq,$ where $p<q$ are two prime numbers satisfying
$q\ne 1(mod\,p).$ If $|G(A)|\ne 1$ then $A$ is a group algebra.
\end{Theorem}

In what follows we find out when $|G(A)|\ne 1$ is guarenteed.
\begin{Proposition}\label{nt}
Let $A$ be a semisimple Hopf algebra 
of dimension $pq,$ where $p<q$ are two prime numbers. If $A^*$ is 
of Frobenius type then either $|G(A)|=p$ and $q=1(mod\,p),$ or
$|G(A)|=pq.$ 
\end{Proposition}
\proof If $A$ is cocommutative then $|G(A)|=pq.$ Otherwise,
$|G(A)|\ne pq,$ and 
by Theorem \ref{notq}, 
$|G(A)|\ne q.$
If $|G(A)|=1,$ then by (\ref{zhu}),
$pq=1+ap^2+bq^2$ for some integers $a,b\ge 0,$ as $A^*$ is of Frobenius 
type. But, $q^2>pq$ hence $b=0$ which yields  
a contradiction.\qed\\

As a corollary we have the following:
\begin{Theorem}\label{frob}
Let $A$ be a semisimple Hopf algebra 
of dimension $pq,$ where $p<q$ are two prime numbers satisfying
$q\ne 1(mod\,p).$ If $A^*$ is of Frobenius type then $A$ is a group algebra.
\end{Theorem}

We consider now semisimple Hopf algebras
$A$ of dimension $pq$ such that $A^*$ is of Frobenius type. 
\begin{Proposition}\label{ca*}
Let $A$ be a non-cocommutative and non-commutative semi-
simple Hopf algebra
of dimension $pq$ over $k,$ 
where $p<q$ are prime numbers. Let $C(A^*)$ be the 
character ring of $A^*.$ 
If $A^*$ is of Frobenius type then:
\ben
\item $C(A^*)$ is commutative.
\item As a coalgebra 
$A=k1\oplus 
kg\oplus \cdots \oplus kg^{p-1}\oplus C_1\oplus\cdots\oplus C_a,$ 
where $a=\frac{q-1}{p}$
and $C_i$ is a simple subcoalgebra of $A$ of dimension $p^2$ for all $1 \le
i \le a.$
\item $gC_i=C_i=C_ig$ for all $1 \le i \le a.$ 
\een
\end{Proposition}
\proof Set $H=kG(A).$ By Theorem \ref{notq} and Proposition \ref{nt}, $dimH=p.$
If $H$ is central in $A$ then $H$ is a 
normal sub Hopf algebra of $A.$ Since $A/AH^+$ is a Hopf algebra of 
dimension $q$ it follows by [Z] that $A/AH^+\cong kC_q,$ and hence 
$A\cong kC_p\#_{\sigma}kC_q.$ Since $kC_p$ 
is central the afforded action is trivial [Ma3, Section 1]. Therefore $A$ is
isomorphic as an algebra to the twisted group ring $kC_p^t[C_q]$ of the
cyclic group $C_q$ over the commutative algebra $kC_p.$ Hence $A$ must be
commutative (See the proof of [Ma2, Proposition 2.3(2)]).
We conclude that $G(A)\cap Z(A)=\{1\}.$

Set $n=dimC(A^*)-1.$ Then, by (\ref{zhu}) there exist two natural 
numbers $a$ and $b$ such that:
\begin{equation} \label{1plus} 
pq=1+ap+bq
\end{equation}
and \begin{equation} \label{nge}
n\ge a+b.
\end{equation}
Clearly, $a\ge 1$ and hence $b<q.$ Moreover, $A^*$ is of Frobenius type 
and $p<q,$ hence by (\ref{square})
\begin{equation} \label{pqge}
pq=p^2(n+1-p)+p.
\end{equation}
Substituting  (\ref{1plus}) and (\ref{nge}) in (\ref{pqge})
yields 
$$pq \ge p^2(a+b+1-p)+p
= p^2\left(\frac{(p-b)q-1}{p}+b+1-p\right)+p$$
and hence $(1-p+b)q\ge (1-p+b)p.$  Since $p<q$ and $b<p,$ this is possible 
if and only if $b=p-1$ and equality holds in (\ref{nge}).
This implies that $C(A^*)$ is commutative and that
$pq=1+ap+(p-1)q$ and hence $a=\frac{q-1}{p}.$
Let 
$$e_{_A},\,e_1,\,\dots,\,e_{_a},\,e_{a+1},\dots ,e_{a+p-1}$$
be the primitive idempotents of $C(A^*),$ where $e_{_A}$ is the
integral of $A$ with $\varepsilon(e_{_A})=1,$ $dim(Ae_i)=p$
for $1\le i\le a$ and $dim(Ae_{a+j})=q$ for $1\le j\le p-1.$
Let $E_0=\frac{1}{p}\sum_{i=0}^{p-1} g^i$ be an integral of $H.$ 
Since $dim(A/AH^+)=q,$ it follows that $AH^+=A(1-E_0)$ has dimension $(p-1)q$ 
and thus $dim(AE_0)=q.$ Moreover, $E_0e_{_A}=e_{_A},$ hence counting
dimensions yields that 
$$E_0=e_{_A}+e_1+\cdots +e_{_a}.$$
Now, $C(A^*)$ is commutative, thus
$dim(C(A^*)e_{_A})=dim(C(A^*)e_i)=1$ for all $1\le i\le a$
and hence 
\begin{equation}\label{eh}
dim(C(A^*)E_0)=a+1.
\end{equation}
Since the set of all the irreducible left $A^*-$modules consists of 
$p$ $1-$dimensional modules and $a=\frac{q-1}{p}$ $p-$dimensional modules 
it follows that
$$A=k1\oplus 
kg\oplus \cdots \oplus kg^{p-1}\oplus C_1\oplus\cdots\oplus C_a$$
as a coalgebra where $C_i$ is a simple
subcoalgebra of $A$ of dimension $p^2$ for all $1\le i\le a.$  
Let 
$$\{1,g,\dots ,g^{p-1},\chi_1,\dots,\chi_a\}$$
be the set of irreducible characters of $A^*,$ where $\chi_i$ corresponds
to $C_i.$ This set clearly forms a basis of $C(A^*).$
Then 
$$C(A^*)E_0=sp\{E_0,E_0\chi_1,\dots,E_0\chi_a\}$$
which implies by (\ref{eh}) that $\{E_0,E_0\chi_1,\dots,E_0\chi_a\}$ 
forms a basis of $C(A^*)E_0.$ 
If $g \chi_i=\chi_j$ for
$i\ne j$ then $E_0\chi_i=E_0\chi_j$ which is a contradiction. Therefore,
$g \chi_i=\chi_i,$ and hence
$$gC_i=C_i=C_ig $$
for all $i.$
This concludes the proof of the proposition. \qed
\begin{Theorem}\label{main}
Let $A$ be a semisimple Hopf algebra
of dimension $pq$ over $k,$  
where $p<q$ are prime numbers.
If $A^*$ is of Frobenius type and $|G(A^*)|\ne 1,$ then $A$ is either a 
group algebra or a dual of a group algebra, and $|G(A)|=p<q$ or $pq.$
\end{Theorem}
\proof  
If $A$ is either cocommutative or commutative then $A$ is either a group  
algebra or a dual of a group algebra respectively. In any event $A^*$ is of
Frobenius type, hence by Proposition \ref{nt}, $|G(A)|=p<q$ or $pq$ and we are 
done.\\
Suppose that $A$ is not cocommutative and not commutative.
Then Proposition \ref{ca*} is applicable. 
Set $H=kG(A).$
By Proposition \ref{nt}, 
$|G(A)|=p.$ By Theorem \ref{notq}, $|G(A^*)|\ne q,$ hence
$|G(A^*)|=p$ too.
Thus we have the following sequence
of maps:
$$H\stackrel{i}{\hookrightarrow}A\stackrel{\pi}{\raro}H$$
where $i$ is the inclusion map and $\pi$ is
a surjection homomorphism of Hopf algebras.
If $\pi\circ i=\varepsilon$ then $H\subseteq K=A^{coH}.$ Since $K$ is a
left coideal of $A,$ it is a direct sum of irreducible left coideals
$K=k1\oplus kg\oplus\cdots\oplus kg^{p-1}\oplus 
V_1\oplus \cdots\oplus V_n,$ where
$dimV_i=p$ for all $i$ as $A^*$ is of Frobenius type. This is a
contradiction since $p$ does not divide
$dimK=q.$
Therefore $\pi\circ i\ne \varepsilon$ and we may assume that $\pi\circ
i=id_H.$ Therefore by Theorem \ref{radf}, there exists $B\subset A$ so
that $A\cong B\times H.$ 
By Theorem \ref{radf}(iv), $B$ is a 
left coideal of $A,$ hence a
direct sum of irreducible left
coideals of $A.$ 
By Proposition \ref{ca*}(2), the dimensions of these irreducible left
coideals are either $1$ or $p.$ 
Since $dimB=q,$ it follows that $B$ contains an irreducible left coideal
$V$ of $A,$ of dimension $p.$
Since $V \subset C$ for some $p^2$-dimensional simple subcoalgebra $C,$ 
it follows by Theorem \ref{radf}(vii) and Proposition \ref{ca*}(3) that 
$V \times H\subseteq C.$ 
But, $dim(V\times H)=p^2=dimC,$ hence $V \times H = C.$
By Theorem \ref{radf}(iii),
$A/AH^+ \cong B$ as coalgebras, and $V$ is the image of $C=V \times H$
under this isomorphism, hence 
$${\rm V\;is\;a\;subcoalgebra\;of\;B.}$$ 

We wish to prove that $V$ is a simple sub
coalgebra of $B$ and thus to reach a contradiction.
Note that since $V$ is an irreducible left coideal of $A$ it follows that
$V\times g^i$ is also an irreducible left coideal of $A$ for all $0\le
i\le p-1.$ By Proposition \ref{ca*}(3), it    
follows that the set 
$$\{V \times g^i|0 \le i \le p-1\}$$ 
is the set of
all the irreducible left coideals of $A$ contained in $C.$ 
Since, $V$ is a
left coideal of $A,$ it 
follows from Theorem \ref{radf}(ii) that $V$ is an $H$ sub-comodule of $B.$ 
Let $\rho:B \raro H \ot B $  
be the comodule structure map, and write
$V=\oplus_{i=0}^{p-1}V_i,$ where $V_i=\rho^{-1} 
(g^i \ot V).$ We claim that $dimV_i=1$ for all
$i.$ Indeed, let $\{v_0,\dots,v_{p-1}\}$ be a basis of $V$ consisting 
of homogenous elements; that is, $\rho(v_i)=g^{m_i}\ot v_i$ for some 
$0\le m_i\le p-1.$ Let $0\ne v \in V$ and write $\gD_B(v)=\sum_{i=0}^{p-1}
b_i \ot v_i.$ Then by Theorem \ref{radf}(vi),
$$\gD_A(v\times 1)=\sum_{i=0}^{p-1} b_i \times g^{m_i}\ot v_i \times 1.$$
Therefore, using Kaplansky's notation [Kap], $L(v \times 1)=
sp\{b_i\times g^{m_i}|0\le i\le p-1\}\subset C$ 
is a {\bf right} coideal
of $A$ of dimension $\le p.$ Since $C$ is 
a simple subcoalgebra of $A$ of dimension $p^2,$ it follows that
$dim(L(v \times 1))=p$ and $L(v \times 1)$ is irreducible.
Therefore by (\ref{inter}), 
$dim(L(v \times 1)\cap(V \times g^i ))=1$ for all $i,$ hence
$\{m_i|0 \le i \le p-1\}=\{0,1,\dots ,p-1\}.$  
Thus $V$ has a basis $\{v_i|0 \le i \le p-1\},$ where $\rho(v_i)=g^i \ot
v_i.$
Since $V$ is an $H$-comodule coalgebra it follows that
$\gD_B(V_i)\subseteq\sum_{j=0}^{p-1}V_j\ot V_{i-j},$ hence 
$\gD_B(v_i)=\sum_{j=0}^{p-1}\ga_{ij}v_j \ot v_{i-j}$ for all $i,$ 
for some $\ga_{ij}\in k.$
Computing $\gD_A(v_i \times 1)$ yields that 
$R_i=L(v_i\times 1)=
sp\{\ga_{ij}v_j\times
g^{i-j}|0 \le j \le p-1\}\subset C$ is a right coideal of $A$ of
dimension$\le p,$ 
for all $i.$ Hence $dimR_i=p$ and
\begin{equation}\label{ri}
R_i=sp\{v_j\times g^{i-j}|0\le j\le p-1\}
\end{equation}
is irreducible. It is straightforward to verify
that $R_i\ne R_t$ for $i\ne t,$ and hence the set $\{R_i|
0 \le i \le p-1\}$ is the set of {\bf all} the irreducible  
right coideals of $A$ which are contained in $C.$

Finally, let $D\subseteq V$
be a subcoalgebra. By Theorem \ref{radf}(vi), $D\times H\subseteq C$ is a
{\bf right} 
coideal of $A$ and hence $D\times H= \oplus_l R_{i_l},$ 
where $R_{i_l}$ is as in (\ref{ri}). But, the image of
$D \times H$ under the map $id \ot \gep:A \raro B$ equals $D,$ while the
image of $\oplus_l R_{i_l}$ under this map equals $V.$ 
Therefore $D=V,$ and hence $V$ is a simple coalgebra. But,
this is a 
contradiction since $dimV=p$ is not a square. 
\qed\\

As a corollary we obtain the following:
\begin{Theorem}\label{rem}
If both $A$ and $A^*$ are of Frobenius type
then $A$ is either a group algebra or the dual of a 
group algebra.
\end{Theorem}
\proof Follows from Proposition \ref{nt} and Theorem \ref{main}.\qed\\

We conclude the paper by considering semisimple Hopf algebras of dimensions
$5p$ and $7p.$
\begin{Lemma}\label{gen}
Let $A$ be a semisimple Hopf algebra of odd dimension over $k.$ If
$|G(A)|=1$ then there exists an irreducible $A^*-$module $V$ with 
$dimV\ge 4.$
\end{Lemma}
\proof
Suppose on the contrary that for any non-trivial $A^*-$irreducible module 
$V,$ $dimV\le 3$. Then by [NR], $dimV=3.$ Hence, $dim(V\ot V^*)=9$ and by 
[La], $V\ot V^*=k\oplus V_1\oplus\cdots\oplus V_i,$ where $V_j\ne k$ is an 
$A^*-$irreducible module for all $j.$ Since $dimV_j=3,$ this is a 
contradiction.\qed
\begin{Theorem}\label{5}
Let $A$ be a semisimple Hopf algebra over $k.$ 
If $dimA=5p,$ $p$ an odd prime, and if
$p=2$ or $4(mod\,5)$ or  
$p\in \{13,23\},$ then $A$ is a group algebra.
\end{Theorem} 
\proof We wish to show that $|G(A)|\ne 1.$
Suppose on the contrary that $|G(A)|=1.$ 
Set $n=dimC(A^*)-1.$ By (\ref{zhu}), there exist two natural
numbers $1\le a$ and $1\le b\le 4$ so that 
\begin{eqnarray*}
& & 5p=1+5a+bp,\\ 
& & n\ge a+b\;\;\; and\\
& & 5p\ge 9(n-1)+16+1
\end{eqnarray*}
where the last inequality follows from (\ref{square}) and Lemma \ref{gen}.
Hence $(-20+9b)p\ge 45b+31.$ But, if $p=2$ or $4(mod\,5)$ then 
$b=2$ or $1$ respectively and if $p\in\{13,23\}$ then $b=3.$ 
In any event this is impossible and we have proved that 
$|G(A)|\ne 1.$  
The result follows now from Theorem \ref{trivial}.\qed
\begin{Theorem}\label{7}
Let $A$ be a semisimple Hopf algebra over $k.$ 
If $dimA=7p,$ $p$ a prime, and if $p=6(mod\,7)$ or $p\in\{17,31\},$
then $A$ is a group algebra. 
\end{Theorem} 
\proof
Suppose $|G(A)|=1$ and set $n=dimC(A^*)-1.$ 
By (\ref{zhu}), there exist two natural
numbers $1\le a$ and $1\le b\le 6$ so that 
\begin{eqnarray*}
& &7p=1+7a+bp,\\
& &n\ge a+b\;\;\;and\\
& &7p\ge 9(n-1)+16+1
\end{eqnarray*}
where the last inequality follows from (\ref{square}) and Lemma \ref{gen}.
Thus, $(-14+9b)p\ge 63b+47.$ But, if $p=6(mod\,7)$ or $p\in\{17,31\}$ 
then this 
is impossible. The result follows now from Theorem \ref{trivial}.\qed
 
\end{document}